\documentclass[10pt]{article}
\usepackage{amsmath,amssymb,amsfonts,amsthm}
\usepackage{epsf,color}
\usepackage{fancyhdr}
\newtheorem{theorem}{Theorem}

\newcommand{\tightoverset}[2]{%
    \mathop{#2}\limits^{\vbox to 0.28ex{\kern -0.2ex\hbox{$#1$}\vss}}}

\usepackage{hyperref}

\vspace{.5cm}

\title{A Note on the Number of Permutations whose Cycle Lengths Are Prime Numbers}

\bigskip

\author{ {\sc Ljuben Mutafchiev}\thanks{American University in
Bulgaria, 1 "Georgi Izmirliev" Square,
 Blagoevgrad 2700, Bulgaria \  E-mail: Ljuben@aubg.edu}\,\,\thanks{Institute of Mathematics and Informatics, Bulgarian Academy of Sciences,
  "Akad. Georgi Bonchev" bl. 8, Sofia 1113, Bulgaria}\,\,\,\thanks{This work was partially supported by Project KP-06-N32/8 with the Bulgarian Ministry of
  Education and Science}}
\date{}
\begin{document}
\maketitle

\begin{abstract}
Let $A$ be a set of natural numbers and let $S_{n,A}$ be the set
of all permutations of $[n]=\{1,2,...,n\}$ with cycle lengths
belonging to $A$. For $A(n)=A\cap [n]$, the limit
$\rho=\lim_{n\to\infty}|A(n)|/n$ (if it exists) is usually called
the density of set $A$. (Here $|B|$ stands for the cardinality of
the set $B$.) Several studies show that the asymptotic behavior of
the cardinality $|S_{n,A}|$, as $n\to\infty$, depends on the
density $\rho$. It turns out that the assumption $\rho>0$ plays an
essential role in the asymptotic analysis of $|S_{n,A}|$. Kolchin
(1999) noticed that there is a lack of studies on classes of
permutations satisfying $\rho=0$ and proposed investigations of
certain particular cases. In this note, we consider the
permutations whose cycle lengths are prime numbers, that is, we
assume that $A=\mathcal{P}$, where $\mathcal{P}$ denotes the set
of all primes. From the Prime Number Theorem it follows that
$\rho=0$ for this class of permutations. We deduce an asymptotic
formula for the summatory function $\sum_{k\le n}
|S_{k,\mathcal{P}}|/k!$ as $n\to\infty$. In our proof we employ
the classical Hardy-Littlewood-Karamata Tauberian theorem.

\end{abstract}

\vspace {.5cm}

 {\bf Key words:} permutation, cycle length, prime number,
 Hardy-Littlewood-Karamata Tauberian theorem
\vspace{.5cm}

{\bf Mathematics Subject classifications:} 05A05, 05A16

\vspace{.2cm}

\section{Introduction, Motivation and Statement of the Main Result}

We start with some notation and conventions that will be used freely
in the text of the paper.

The letter $p$ without subscript denotes a prime number. We write
$\mathcal{P}$ for the set of all primes. If the primes in
$\mathcal{P}$ are arranged in increasing order, then $p_k$ denotes
the $k$th smallest prime, $k=1,2,...$. For a function $\psi(y)$ of
a real variable $y$ which assumes real or complex values,  we
write $\sum_p \psi(p)$ instead of $\sum_{p\in\mathcal{P}}
\psi(p)$. We also use the usual notation $\pi(y)$ for the number
of primes not exceeding $y>0$.

Further on, $A$ always denotes a set of positive integers. Let
$[n]=\{1,2,...,n\}$. Then, we define the set $A(n)=A\cap [n]$. We
denote the cardinality of any set $B$ by $|B|$. The limit
\begin{equation}\label{density}
\rho=\lim_{n\to\infty}\frac{|A(n)|}{n},
\end{equation}
if it exists, is called the density of the set $A$. Clearly,
$0\le\rho\le 1$.

We write $S_{n,A}$ for the set of permutations of $n$ letters whose
cycle lengths belong to $A$. We denote the cardinality of $S_{n,A}$
by $P_{n,A}=|S_{n,A}|$. If $A=\mathcal{P}$, we write
$P_n=|S_{n,\mathcal{P}}|$.

For two sequences of real numbers $\{u_n\}_{n\ge 1}$ and
$\{v_n\}_{n\ge 1}$, we write $u_n\sim v_n$ if $\lim_{n\to\infty}
u_n/v_n=1$. In a similar and traditional manner, we shall also
explore the commonly used symbol $O(.)$; for its definition, we
refer the reader, e.g., to \cite{D58}, \cite{O95} and \cite{T15}.

Finally, we show the approximations of two constants that will be
used further, namely: the Euler's constant
$$
\gamma=0.577215
$$
and the constant
\begin{equation}\label{tenen}
c=\gamma-\sum_p\left(\log{\left(\frac{1}{1-1/p}-\frac{1}{p}\right)}\right)
=0.261497,
\end{equation}
which appears in an asymptotic estimate of the partial sum
$\sum_{p\le y} 1/p$ and in Mertens' formula for the partial
product $\prod_{p\le y}(1-1/p)$ as $y\to\infty$. More details may
be found in \cite[Theorems 1.10, 1.12]{T15}.

The problem on the asymptotic enumeration of permutations whose
cycle lengths are constrained has attracted the interest of
several authors in the middle of the last century. The last two
chapters of Kolchin's monogrph \cite{K99} are devoted to certain
aspects from this topic. Kolchin surveyed various important
results and discussed in detail the existing approaches and
methods applied in this area. It seems that one of the main
reasons for the interest in asymptotic enumeration problems of
this type is the relationship between the cardinalities of sets of
permutations $S_{n,A}$, for certain particular sets $A$, and the
theory of equations containing an unknown permutation of $n$
letters. For more details on asymptotic enumeration problems
related to $S_{n,A}$, we refer the reader, e.g., to \cite[Section
4.4 and Chater 5]{K99}, \cite[Section 8.2]{O95} and \cite[Chapter
3]{Y05}.

The enumeration of permutations with cyclic structure constrained
by a certain set $A$ is based on the following generating
function:
\begin{equation}\label{arbgf}
f_A(z):=\sum_{n=0}^\infty\frac{P_{n,A}z^n}{n!}
=\exp{\left(\sum_{k\in A}\frac{z^k}{k}\right)}.
\end{equation}
A proof of this identity may be found, e.g., in \cite[Theorem
5.1.2]{K99}. It turns out that, for several $A$'s, the asymptotic
behavior of $P_{n,A}$ depends on the density $\rho$, given by
(\ref{density}). The case $\rho>0$ was studied under several
additional conditions on $A$ by many authors. As an illustration,
below we give a typical asymptotic result obtained by Yakimiv  in
\cite{Y91} (see also \cite[Theorem 3.3.1]{Y05}).

{\bf Theorem (Yakymiv (1991)).} {\it Suppose that, for a certain
$A$, its density is $\rho>0$. Moreover, for $m\ge n$ and $m=O(n)$,
we assume that
\begin{equation}\label{sq}
|\{k: k\le m, k\in A, m-k\in A\}|/n\to\rho^2
\end{equation}
as $n\to\infty$. Then, we have
\begin{equation}\label{yakasymp}
P_{n,A}\sim n! n^{\rho-1}e^{L(n)-\gamma\rho}/\Gamma(\rho),
\end{equation}
where $L(n)=\sum_{k\in A(n)} 1/k-\rho\log{n}$ and $\Gamma(.)$ is the
Euler's gamma function.}

{\it Remark 1.} Condition (\ref{sq}) shows that the set $A$ may be
considered as a realization of a random set containing each
integer with probability $\rho$ independently from the other
integers. It is also interesting to note that, under entirely
different conditions on $A$, Kolchin \cite{K92} (see also
\cite[Theorem 4.4.10]{K99}) obtained the same asymptotic
equivalence (\ref{yakasymp}).

There are many examples of permutations with constraints on their
cycle structure whose basic generating function (\ref{arbgf}) has
a single dominant singularity at $z=1$. In Sections 11 and 12 of
his excellent survey \cite{O95}, Odlyzko classified and discussed
the methods used to extract asymptotic information about
coefficients of analytic generating functions. He suggested two
main classes of generating functions, depending on whether the
main singularity is large or small. Functions with large
singularities (i.e., ones that grow rapidly as the argument
approaches the circle of convergence) are usually analyzed using
the saddle point method; see \cite[Sections 12.1 and 12.2]{O95}.
In \cite[Chapters 4 and 5]{K99} Kolchin demonstrated this method
several times in the context of permutations with constraints on
their cycles. On the other hand, some generating functions of the
form (\ref{arbgf}) have small singularities on the circle of
convergence and admit applications of other methods: Tauberian
theory and transfer theorems (including those due to Darboux and
Jungen); see \cite[Sections 8.2, 11.1 and 11.2]{O95} and
\cite[Chapter 3]{Y05}.

In \cite[Section 4.4]{K99}, Kolchin surveyed asymptotic results on
constrained permutations whose density $\rho$ is positive and
noticed that there is a lack of studies in the case $\rho=0$.
Kolchin concluded his discussion with an open problem: he proposed
to handle particular cases of $A$'s whose density $\rho=0$. This,
in fact, motivates us to study the case $A=\mathcal{P}$ and
$P_{n,A}=P_n$. In this case the Prime Number Theorem (PNT) implies
that $\rho=0$. (For more details, historical remarks and stronger
results on the remainder term in the PNT, we refer the reader to
\cite[p. 12 and Chapter II.4]{T15}.) We also note that $P_n/n!$ is
the probability that a permutation of $n$ letters chosen uniformly
at random has cycle lengths which are prime numbers. The sequence
$\{P_n\}_{n\ge 1}$ is A218002 in the On-Line Encyclopedia of
Integer Sequences \cite{S}.

Below we state our main result.

\begin{theorem} We have
$$
\sum_{k=0}^n \frac{P_k}{k!}\sim e^c\log{n}
$$
as $n\to\infty$, where $c$ is the constant, given by (\ref{tenen})
and $e^c=1.298873$.
\end{theorem}

{\it Remark 2.} Theorem 1 has a probabilistic interpretation.
Suppose that the set $[k], k=1,2,...,n,$ is chosen uniformly at
random (i.e., with probability $1/n$) and then a permutation of
$[k]$ is selected uniformly at random (i.e. with probability
$1/k!$). Combining Theorem 1 with the total probability formula,
one can observe that the probability that this permutation has
cycle lengths, which are prime numbers, is asymptotically
equivalent to
$$
\frac{e^c\log{n}}{n}.
$$

The proof of Theorem 1 relies on a classical Tauberian theorem of
Hardy, Littlewood and Karamata (see \cite{H49} and \cite[Theorem
8.7]{O95}).

Our paper is organized as follows. The next Section 2 contains
some remarks on our method of proof. We discuss there the
singularities of the underlying generating function and answer the
question: why we have chosen Tauberian theory as a main tool of
our approach? The proof of Theorem 1 is given in Section 3.

\section{Some Remarks on the Method of Proof}

We first note that, for $A=\mathcal{P}$ and $P_{n,A}=P_n$,
(\ref{arbgf}) implies
\begin{equation}\label{gfprime}
f(z):=\sum_{n=0}^\infty\frac{P_n z^n}{n!}
=\exp{\left(\sum_p\frac{z^p}{p}\right)}.
\end{equation}
For the sake of convenience, we also set
\begin{equation}\label{phi}
\varphi(z):=\sum_p\frac{z^p}{p}=\sum_{n=1}^\infty\frac{z^{p_n}}{p_n},
\end{equation}
where in the last equality we have used that the primes in
$\mathcal{P}$ are arranged in increasing order. As first
observation, we will show that the unit circle is a natural
boundary for $\varphi(z)$. The proof of this fact relies
essentially on the Fabry gap theorem, which we state below. For
more details, we refer the reader to \cite[Section 5.3]{M94}.

{\bf Fabry Gap Theorem.} {\it Suppose that the power series
$$
g(z)=\sum_{n=1}^\infty a_n z^{\lambda_n}
$$
has radius of convergence $1$, where
$0\le\lambda_1<\lambda_2<...$. If
$lim_{n\to\infty}\lambda_n/n=\infty$, then the circle of
convergence $|z|=1$ is a natural boundary of $g$.}

It is well known that from the PNT in its minimal form, namely
$\pi(y)\sim y/\log{y}$ as $y\to\infty$, it follows that the $n$th
smallest prime $p_n$ satisfies
\begin{equation}\label{pn}
p_n\sim n\log{n} \quad n\to\infty;
\end{equation}
see, e.g., \cite[Part I, Chapter 1.0, Exercise 14]{T15}. The gaps
between the successive primes in $\mathcal{P}$ and (\ref{pn}) show
that $\varphi$ satisfies the conditions of the Fabry gap theorem.
Thus we observe that the circle $|z|=1$ is a natural boundary for
$\varphi$ and $f$. Hence any application of transfer theorems (see
\cite[Section 11.1]{O95}) fails since an analytic continuation of
the underlying generating function beyond its circle of
convergence is not possible. We recall that the general transfer
theorems (see \cite[Theorem 11.4]{O95}) require that the function
$f(z)$, defined by (\ref{gfprime}), has to be analytic in the
domain $\Delta(1,\phi,\eta)=\{z:|z|\le
1+\eta,|\arg{(z-1)}|\ge\phi\}\setminus \{1\}$, where $\eta>0$ and
$0<\phi<\pi/2$. The dominant singularity of the function $f(z)$
does not satisfy also the algebraic-logarithmic conditions of the
Darboux and Jungen theorems \cite[Section 11.2]{O95} and thus
their application is also impossible. Both methods discussed above
are usually used to extract asymptotic information for the
coefficients of analytic generating functions with small dominant
singularities. As it was expected, our attempt to imitate a saddle
point approximation failed since the function $f(z)$ does not have
large enough singularity at $z=1$. The computational details for
showing this are based on the PNT. In fact, in \cite{M15} was
shown that a weak form of the PNT (namely,
$\pi(y)=y/\log{y}+O(y/\log^2{y}), y\to\infty$) yields the
asymptotic of the first derivative of the function $\varphi(z)$
(see (\ref{phi})), as $z\to 1^-$. We have
$$
\varphi^\prime(z)\sim\frac{1}{(1-z)\log{\frac{1}{1-z}}}.
$$
(In addition, De Bruijn \cite[Section 3.14, Exercise 5]{D58}
proposed an asymptotic series in the sense of Poincar\'{e} for
$\varphi^\prime(z)$.) In a way similar to that given in
\cite{M15}, one can also obtain
$$
\varphi^{\prime\prime}(z)\sim\frac{1}{(1-z)^2\log{\frac{1}{1-z}}}
$$
and
$$
\varphi^{\prime\prime\prime}(z)\sim\frac{2}{(1-z)^3\log{\frac{1}{1-z}}}.
$$
These computations show that the third-order term in the Taylor
expansion of $\varphi(Re^{i\theta})$ around the point $\theta=0$
is too large as $R\to 1^-$, and so, we are in the situation of an
invalid application of the saddle point method, described in
\cite[Example 12.3]{O95}.

The disadvantages discussed above show that these methods are not
adequate for asymptotic analysis of the coefficients of the function
$f(z)$, given by (\ref{gfprime}). Clearly, its singularity at $z=1$
is small. Hence we prefer to apply a classical Tauberian theorem by
Hardy-Littlewood-Karamata. Its proof may be found in \cite[Chapter
7]{H49}. We use it in the form given in \cite[Section 8.2]{O95}.

{\bf Hardy-Littlewood-Karamata Theorem.} {\it (See [7; Theorem
8.7, p. 1225].) Suppose that $h_k\ge 0$ for all $k$, and that
$$
h(z)=\sum_{k=0}^\infty h_k z^k
$$
converges for $0\le z<r$. If there is an $\alpha\ge 0$ and a
function $L(t)$ that varies slowly at infinity such that
\begin{equation}\label{funcsim}
h(z)\sim (r-z)^{-\alpha}L\left(\frac{1}{r-z}\right), \quad z\to
r^-,
\end{equation}
then
\begin{equation}\label{sumsim}
\sum_{k=0}^n h_k r^k\sim\left(\frac{n}{r}\right)^\alpha
\frac{L(n)}{\Gamma(\alpha+1)}, \quad n\to\infty.
\end{equation}}

{\it Remark 3.} A function $L(t)$ varies slowly at infinity if,
for every $u>0$, $L(ut)\sim L(t)$ as $t\to\infty$.

{\it Remark 4.} Odlyzko \cite[Example 8.8]{O95} applied the
Hardy-Littlewood-Karmata Tauberian theorem to obtain a general
formula for the $n$th partial sum of the coefficients
(probabilities) $P_{k,A}/k!$, where the set $A$ has density
$\rho\in [0,1]$. He showed that
$$
\sum_{k=0}^n\frac{P_{k,A}}{k!}\sim\frac{f_A(1-1/n)}{\Gamma(\rho+1)},
$$
where $f_A$ is the function, defined by (\ref{arbgf}). Hence, an
essential step in the proof of Theorem 1 is to determine the
asymptotic behavior of $f(z)$ as $z\to 1^-$ (see (\ref{gfprime})).
Additional results and references on studies in this direction may
be also found in \cite{P92}.

\section{Proof of Theorem 1}

For convenience, in (\ref{gfprime}) and (\ref{phi}) we set
$z=e^{-t}$. Then, we deal first with $\varphi(e^{-t})$. We need to
study the behavior of $\varphi(e^{-t})$ as $t\to 0^+$, since by
(\ref{phi}) $z\to 1^-$ and $t\to 0^+$ are equivalent. We start our
computations with the following representation:
\begin{equation}\label{rep}
\varphi(e^{-t}) =\sum_p\frac{1-(1-e^{-pt})}{p} =\varphi_1(e^{-t})
+\varphi_2(e^{-t})+\varphi_3(e^{-t}),
\end{equation}
where

$$
\varphi_1(e^{-t})=
\sum_{p\le\left(\frac{1}{t}\log{\frac{1}{t}}\right)/\log{\log{\frac{1}{t}}}}
\frac{1}{p},
$$

$$
\varphi_2(e^{-t})=
-\sum_{p\le\left(\frac{1}{t}\log{\frac{1}{t}}\right)/\log{\log{\frac{1}{t}}}}
\frac{1-e^{-pt}}{p},
$$

$$
\varphi_3(e^{-t})=
\sum_{p>\left(\frac{1}{t}\log{\frac{1}{t}}\right)/\log{\log{\frac{1}{t}}}}
\frac{e^{-pt}}{p}.
$$

We notice that $\varphi_j(e^{-t}), j=1,2,3$, are well defined for
sufficiently small $t$. (An upper bound for $t$ could be $e^{-e}$
since, for $t<e^{-e}$, the function
$\log{\log{\log{\frac{1}{t}}}}$, which appears in our further
computations, is well defined.)

The asymptotic estimate for $\varphi_1(e^{-t})$ is obtained by
means of a classical result from Number Theory, given in
\cite[Theorems 1.10 and 1.12]{T15} and related to the asymptotic
growth of the partial sum $\sum_{p\le y} 1/p$ as $y\to\infty$. An
application of this result immediately yields
\begin{eqnarray}\label{one}
& & \varphi_1(e^{-t}) =
\log{\log{\left(\frac{\log{\frac{1}{t}}}{t\log{\log{\frac{1}{t}}}}\right)}}
+c +O\left(\frac{1}{\log{\frac{1}{t}}}\right) \nonumber \\
& & =\log{\left(\left(\log{\frac{1}{t}}\right)\left(1
+\frac{\log{\log{\frac{1}{t}}}-\log{\log{\log{\frac{1}{t}}}}}{\log{\frac{1}{t}}}\right)\right)}
+c +O\left(\frac{1}{\log{\frac{1}{t}}}\right) \nonumber \\
& & =\log{\log{\frac{1}{t}}}
+\log{\left(1+\frac{\log{\log{\frac{1}{t}}}-\log{\log{\log{\frac{1}{t}}}}}{\log{\frac{1}{t}}}\right)}
+c +O\left(\frac{1}{\log{\frac{1}{t}}}\right) \nonumber \\
& & =\log{\log{\frac{1}{t}}}
+O\left(\frac{\log{\log{\frac{1}{t}}}}{\log{\frac{1}{t}}}\right)
+c +O\left(\frac{1}{\log{\frac{1}{t}}}\right) \nonumber \\
& & =\log{\log{\frac{1}{t}}}+c+
O\left(\frac{\log{\log{\frac{1}{t}}}}{\log{\frac{1}{t}}}\right),
\end{eqnarray}
where the constant $c$ is defined by (\ref{tenen}).

The asymptotic estimate for $\varphi_2(e^{-t})$ will be based
again on the PNT in its minimal form. We apply first the
inequality $(1-e^{-pt})/p\le t$ to each term in
$\varphi_2(e^{-t})$. Then, by the PNT, we have
 \begin{equation}\label{two}
|\varphi_2(e^{-t})|\le
t\sum_{p\le\left(\frac{1}{t}\log{\frac{1}{t}}\right)/\log{\log{\frac{1}{t}}}}
1 \sim\frac{t\left(\frac{1}{t}\log{\frac{1}{t}}\right)}
{\left(\log{\log{\frac{1}{t}}}\right)\log{\frac{1}{t}}}
=O\left(\frac{1}{\log{\log{\frac{1}{t}}}}\right)
\end{equation}
as $t\to 0^+$. Finally, for $\varphi_3(e^{-t})$, we obtain
\begin{eqnarray}
& & \varphi_3(e^{-t})\le \frac{t\log{\log{\frac{1}{t}}}}{\log{1/t}}
\sum_{p>\left(\frac{1}{t}\log{\frac{1}{t}}\right)/\log{\log{\frac{1}{t}}}}
e^{-pt} \nonumber \\
& & <\frac{t\log{\log{\frac{1}{t}}}}{\log{{1/t}}}
\sum_{m>\left(\frac{1}{t}\log{\frac{1}{t}}\right)/\log{\log{\frac{1}{t}}}}
e^{-mt}, \nonumber
\end{eqnarray}
where in the last sum we have also included all terms $e^{-mt},
m\notin\mathcal{P}$. The last sum is a tail of a geometric
progression and hence we easily obtain the following asymptotic
estimate:
\begin{equation}\label{three}
\varphi_3(e^{-t})
=O\left(\frac{\log{\log{\frac{1}{t}}}}{\log{1/t}}
e^{-\left(\log{\frac{1}{t}}/\log{\log{\frac{1}{t}}}\right)}\right).
\end{equation}
Combining (\ref{rep}) - (\ref{three}) all together, we deduce that
$$
\varphi(e^{-t})=\log{\log{\frac{1}{t}}}+c
+O\left(\frac{1}{\log{\log{\frac{1}{t}}}}\right), \quad t\to 0^+.
$$
Then (\ref{gfprime}) and (\ref{phi}) imply
\begin{equation}\label{ft}
 f(e^{-t}) =e^c \left(\log{\frac{1}{t}}\right)
\left(1+O\left(\frac{1}{\log{\log{\frac{1}{t}}}}\right)\right),
\quad t\to 0^+.
\end{equation}
To complete the proof of Theorem 1, we recall that
$t=\log{\frac{1}{z}}$. Since
$$
\log{\frac{1}{z}}=-\log{z}=-\log{(1-(1-z))}\sim 1-z, \quad z\to
1^-,
$$
we can rewrite (\ref{ft}) as follows
$$
f(z)\sim e^c\log{\frac{1}{1-z}}, \quad z\to 1^-.
$$
Therefore the series $f(z)$, defined by (\ref{gfprime}), satisfies
condition (\ref{funcsim}) of Hardy-Littlewood-Karamata Tauberian
theorem with $h(z)=f(z), r=1, \alpha=0$ and $L(t)=\log{t}$. Now, the
asymptotic equivalence of Theorem 1 follows immediately from
(\ref{sumsim}).

\end{document}